\documentclass[letterpaper,conference,10pt]{ieeeconf}
\IEEEoverridecommandlockouts
\pdfoutput=1

\makeatletter
\let\NAT@parse\undefined
\makeatother

\usepackage[square, compress, numbers]{natbib}

\usepackage{amsmath,graphicx,epsfig,color,amsfonts, algorithm, subfigure,latexsym}

\graphicspath{{./Figs/}}

\usepackage[noend]{algorithmic}

\def\real{\mathbb{R}}

\newcommand{\setdef}[2]{\{#1 \; | \; #2\}}

\newcommand\oprocendsymbol{\hbox{$\square$}}
\newcommand\oprocend{\relax\ifmmode\else\unskip\hfill\fi\oprocendsymbol}

\renewcommand{\theenumi}{\roman{enumi}}

\newcommand\bit[1]{\textit{\textbf{#1}}}

\def \bs {\boldsymbol}
\def \mc {\mathcal}

\def \etal {\emph{et al.}}

\newtheorem{theorem}{Theorem}
\newtheorem{proposition}[theorem]{Proposition}

\newtheorem{example}{Example}

\newtheorem{claim}{Claim}

\renewcommand{\theenumi}{(\roman{enumi}}

\title{A Realization Theory for Bio-inspired Collective Decision-Making \thanks{This research has been supported in part by ONR grant 
N00014-14-1-0635, ARO grant W911NF-14-1-0431 and NSF grant ECCS-1135724.}}

\author{Alessio~Franci~\hspace{0.5in}~Vaibhav~Srivastava~\hspace{0.5in}~Naomi Ehrich Leonard
\thanks{A. Franci is with the Department of Engineering, University of Cambridge, Cambridge, UK,  \texttt{af529@cam.ac.uk}. }
\thanks{V. Srivastava and N. E. Leonard are with the Department of Mechanical and Aerospace Engineering, Princeton University, Princeton, NJ, USA, \texttt{\{vaibhavs, naomi\}@princeton.edu}.
}
}

\begin{document}

\maketitle

\begin{abstract}

The collective decision-making exhibited by animal groups provides enormous inspiration for multi-agent control system design as it  embodies several features that are desirable in engineered networks, including robustness and adaptability, low computational effort, and an intrinsically decentralized architecture. However, many of the mechanistic models for collective decision-making are described at the population-level abstraction and are challenging to implement in an engineered system. We develop simple and easy-to-implement models of opinion dynamics that realize the empirically observed collective decision-making behavior as well as the behavior predicted by existing models of animal groups. Using methods from Lyapunov analysis, singularity theory, and monotone dynamical systems, we rigorously investigate the steady-state decision-making behavior of our models. 
 
\end{abstract}

\section{Introduction}

Collective decision-making in animal groups has received significant attention in a broad scientific community~\cite{LC-CL:09, IDC:09}.    Animals rely on their ability to make decisions among alternatives; many species perform remarkably well, even in the presence of disturbances, changes in the environment, and limitations on individual sensing, communicating, and computing.  For example, fish schools make good collective choices among directions to forage \cite{IDC-JK-NRF-SAL:05}, migrating bird flocks  choose efficiently when to depart  as a group from a rest stop \cite{Eikenaar2014}, and honeybee swarms jointly select the best option for a new nest \cite{Seeley12}.  To explain the remarkable group-level behaviors, it has been hypothesized that animals behave like a networked multi-agent system, using decentralized strategies and interactions based on social cues.      


Models for collective decision-making 
in scenarios such as honeybee nest site selection~\cite{Seeley12, Pais2013} have been developed that exhibit many desirable features including the unanimity of decisions, the robustness and adaptability to variations in the environment,  and the sensitivity of the decision-making process to the absolute values, and the differences between  values, of the alternatives.  However, many of these  are population-level models and  do not provide insights on the role of the network structure on the decision-making process.  Nor are they easy to  implement in an engineered network.

Mechanistic decision-making models can be used to study the decision-making process. In deterministic mechanistic models, bifurcations are the most fundamental phenomena that capture the sensitivity of the decision-making to  environmental factors.

In this paper, we develop  simple agent-level models that realize the collective decision-making behavior predicted by  population-level models as well as empirical data and agent-based models of animal groups in two alternative choice tasks.   Our proposed model provides a framework for systematic bio-inspired design of collective decision-making.  It can be used for implementing high performing control strategies in an engineered network, and for systematic model-based investigation of collective decision-making of animal groups, which may lead to new testable hypotheses. 



Couzin~\etal~\cite{IDC-JK-NRF-SAL:05} developed a mechanistic model to examine the role of  preferences in a moving animal group choosing between two alternatives.  Numerical simulations suggested a bifurcation in behavior from a  decision-deadlock to a group decision; this was studied analytically using a deterministic model in~\cite{BN-NEL-etal:09, NEL-TS-etal:12}.   In \cite{Couzin2011}, schooling fish data were analyzed from experiments with an asymmetry in preference strength and manipulation of the fraction of the group with no preference.  The model fitted to the data showed a perturbation (or {\it unfolding} \cite{Golubitsky1985}) of the pitchfork bifurcation corresponding to a decision-deadlock  for small fraction with no preference, and a decision for one option for high fraction with no preference.

For the honeybee nest site selection problem between two nest sites, Seeley~\etal~\cite{Seeley12} developed a mechanistic model in which the decision-making process is captured through a pitchfork bifurcation. Their model involves a cross-inhibitory signal that serves as the bifurcation parameter. For small values of the cross-inhibitory signal, there is a decision-deadlock, while for large values of the signal, a unanimous decision is made for one of the nest sites. If the value of each nest site is the same, then a symmetric pitchfork bifurcation is observed, otherwise an unfolding of the pitchfork bifurcation is observed. Pais~\etal~\cite{Pais2013} performed a rigorous analysis of this model and drew connections to behavioral models of decision-making proposed in~\cite{RB-EB-etal:06}. 

Models for opinion dynamics and information assimilation in social networks~\citep{JL:07, VDB-JMH-JNT:09, AM-FB:12,AN-BT:12, DA-AO-AP:10, AJ-AS-ATS:10} are receiving growing attention in the control community. These models are fairly general and exhibit very rich behavior. However, they do not incorporate the decision-making process nor capture the empirically observed collective decision-making behavior in animal groups. 
Behavioral models for decision-making in two-alternative-choice tasks have been integrated with opinion dynamics models to capture collective decision-making in groups~\cite{VS-NEL:13f}. However, generic mechanistic agent-based models for collective decision-making have not been much explored. 

The models in this paper comprise interconnections of linear time invariant systems and sigmoidal nonlinearities.
It has been shown in~\cite{AF-RS:14} that such systems exhibit rich nonlinear dynamics, yet, are analytically tractable and easily implementable in real world systems.

The major contributions of this paper are threefold. First, we propose and analyze an uninformed opinion dynamics model. We show that for small values of social effort, the stable fixed points of the model correspond to decision-deadlock; while for values of social effort above a critical value, the model admits a pitchfork bifurcation and  a unanimous decision is achieved. Our analysis relies on Lyapunov theory, singularity theory, and monotone system theory. 

Second, we extend the uninformed opinion dynamics to the informed opinion dynamics in which some agents have access to the external stimuli. Using Lyapunov-Schmidt reduction, we show that the informed opinion dynamics are an unfolding of the pitchfork bifurcation obtained in the uninformed case. 
Third, we apply the informed opinion dynamics model to the nest site selection problem in honeybees and show that our model captures the predictions of the associated population-level models. 

The remainder of the paper is organized as follows. We present our models for uninformed and informed opinion dynamics in \S\ref{SEC: model derivation}. We review fundamentals of singularity theory and Lyapunov-Schmidt reduction in~\S\ref{sec:review}. Our analysis of uninformed and informed opinion dynamics is contained in \S\ref{SEC: gloabl local analysis}~and~\S\ref{sec:unfold}, respectively. We apply the informed opinion dynamics model to the honeybee nest site selection problem in~\S\ref{SEC: honey bee}. Finally, we conclude in \S\ref{SEC: conclusions}.

\section{A realization of collective decision-making}
\label{SEC: model derivation}

We consider a set of $N$ interconnected agents. Let $A\in\mathbb R^{N\times N}$ be the agent network adjacency matrix, with $a_{ij}\geq0$ and $a_{ii}=0$ for all $i,j=1,\ldots,N$ and $j\neq i$. Let $D\in\mathbb R^{N\times N}$ be a  diagonal matrix with $D_{ii}=d_i:=\sum_{j=1}^N a_{ij}$ and let $L=D-A$ be the network Laplacian matrix. We make the {\it standing assumption} that the interconnection graph is strongly connected and balanced, that is, $\text{rank}\ L=N-1$ and $\bs 1_N^TL=L \bs 1_N=0$, where $\bs 1_N$ is the $N$-column-vector with all unitary entries.

In this paper, we study the following distributed dynamics as a candidate  for the realization of bio-inspired collective decision-making behavior
\begin{equation}\label{EQ: LTIs network}
\dot x_i=-d_ix_i+\sum_{j=1}^N ua_{ij} S(x_j),\quad i=1,\ldots,N
\end{equation}
where $S:\mathbb R\to (-1, 1)$ is a smooth sigmoidal function that satisfies the following conditions:
\begin{enumerate}
\item  $S(-z)=-S(z)$, for all $z\in\mathbb R$ (odd);
\item $S^{(1)}(z)>0$, for all $z\in\mathbb R$ (monotone);
\item  $\lim_{z\to\pm\infty}S^{(1)}(z)=0$\ref{EQ: LTIs network} (saturated);
\item For all $n\in\mathbb N$, $S^{(2n+1)}(0)\neq 0$ and $S^{(1)}(0)=1$;
\item For all $n\in\mathbb N$ and $z \ne 0$, ${\rm sgn}\left(S^{(2n)}(z)\right)=-{\rm sgn}(z)$. 
\end{enumerate}
 In  vector form, dynamics~\eqref{EQ: LTIs network} can be written as
\begin{equation}\label{eq:network-ltis-vector}
\dot{\bs x} = - D \bs x + u A \bs S(\bs x),
\end{equation}
where $\bs S(\bs x)$ is a $N$-column-vector with $i$-th entry $S(x_i)$.  We interpret the state $\bs x$ as the vector of agent opinions and treat it as a measure of ensuing decisions. Accordingly,   we refer to~\eqref{eq:network-ltis-vector} as the \emph{uninformed opinion dynamics}.

The uninformed opinion dynamics~(\ref{EQ: LTIs network}) is a network of output-saturated first order LTI systems\footnote{Many functional forms are possible for which all the results of the paper remains true (input saturation, etc.). We pick (\ref{EQ: LTIs network}) to keep the illustration concrete.}. This model may be interpreted in the following way.
The term $u S(x_j)$ in~\eqref{EQ: LTIs network} is the opinion of agent $j$ as perceived by a generic agent $i$. In particular, 
the saturation function models an agent's assessment of the opinion of other agents: opinions with small values are assessed as they are; while the opinions with large values  are assessed to have a smaller value. The parameter $u$ controls this smaller value and models the social effort: higher social effort leads to a broader range of opinions being assessed correctly. 
With this interpretation, the uninformed opinion dynamics~\eqref{EQ: LTIs network} is the continuous time version of the process in which each agent at each time updates her opinion to a convex combination of her opinion with the perceived opinions of her neighbors.

Heuristically, the choice of model (\ref{EQ: LTIs network}) is dictated by the following observation. The linearization of model (\ref{EQ: LTIs network}) at the origin for $u=1$ is a consensus dynamics (with Laplacian $L=A-D$) and, therefore, it is singular with a single zero eigenvalue along the consensus manifold. This implies the presence of a bifurcation with center space tangent to the consensus manifold  \cite[Theorem 3.2.1]{Guckenheimer2002}, which will generically be a pitchfork. This claim and the fact that this bifurcation globally determines (or ``organizes") the dynamics of (\ref{EQ: LTIs network}) are proved in Theorem \ref{thm: main unp pitch} in Section \ref{SEC: gloabl local analysis}.

In the following we also study a version of~\eqref{eq:network-ltis-vector} in which an external stimuli is added. This external stimuli may represent agent preferences or their exposure to some external information source.
 We refer to such dynamics as the \emph{informed opinion dynamics} and define it by 
\begin{equation*}
\dot x_i=-d_ix_i+\sum_{j=1}^Nua_{ij}S(x_j)+\alpha_i,\quad i=1,\ldots,N.
\end{equation*}
In  vector form the informed opinion dynamics are 
\begin{equation}\label{eq:opinion-forced}
\dot{\bs x} = - D \bs x + u A \bs S(\bs x) + \bs \alpha,
\end{equation}
where $\bs \alpha \in \real^N$ is the preference term. 

The informed opinion dynamics model~\eqref{eq:opinion-forced}  can be interpreted similarly to the uninformed opinion dynamics model. In particular, the uninformed opinion dynamics~\eqref{eq:opinion-forced} is the continuous time version of the process in which each agent at each time (i) computes a convex combination of her opinion with the perceived opinions of her neighbors, and (ii) updates her opinion by the sum of the convex combination and the external stimuli. 


\section{A review of singularity theory and Lyapunov-Schmidt reduction} \label{sec:review}

In this section we rapidly review some of the fundamentals of singularity theory and Lyapunov-Schmidt reduction. The reader is referred to~\cite{Golubitsky1985} for a comprehensive treatise on these topics. See also \cite{AF-RS:14} for a short introduction with a control theoretical flavour.

\subsection{Singularity and universal unfolding}

A scalar bifurcation problem
\[
g(y, u)=0,
\]
is defined by the zero set of the smooth function $g$ as the \emph{bifurcation parameter} $u \in \real$ is varied. The independent variable $u$ usually models a control parameter, whereas the dependent variable $y\in\real$ is a measured quantity of interest. The set $\{(y,u)\ |\ g(y,u)=0 \}$ is called the \emph{bifurcation diagram}. A zero $(y^*, u^*)$ is said to be a \emph{bifurcation point} if the bifurcation diagram is not regular at $(y^*, u^*)$. For instance, the number of zeros might change at a bifurcation point. A necessary condition for $(y^*, u^*)$ to be a bifurcation point is $g_y(y^*, u^*)=0$, where $g_y$ denotes the derivative with respect to $y$. Such a point is called a \emph{singularity}.  If at $(y, u) = (y^*, u^*)$,
\[
g = g_y = g_{yy} = g_{u} =0, \; g_{yyy}>0, \text{ and } g_{u y} <0,
\]
then the singularity is a \emph{pitchfork}. The above conditions are said to solve the \emph{recognition problem} of the pitchfork. At a pitchfork singularity, the bifurcation diagram exhibits the typical pitchfork shape ( $-\hspace{-2mm}\in$ ), corresponding to a single zero bifurcating into three zeros.

Let $G(y, u, \bs \alpha)$, $\bs \alpha \in \real^k$ be a $k$-parameter family of bifurcation problems such that $G(y, u, 0) = g(y, u)$; and, given any smooth perturbation term $\epsilon p (y, u,  \bs \alpha)$ with $\epsilon$ sufficiently small, there exists $\bs \alpha$ such that $g + \epsilon p$ and $G(y, u, \bs \alpha)$ are {\it strongly equivalent}, that is, $G$ is obtained from $g + \epsilon p $ via a local diffeomorphism of the form $(y, u)\mapsto(Y(y,u),U(u))$ and multiplication by a positive function. Then, $G$ is said to be a \emph{universal unfolding} of $g$ if $k$ is equal to the {\it codimension} of $g$. The family of bifurcation diagrams associated with a universal unfolding are $\setdef{ (y, u) }{G(y, u, \bs \alpha) =0,\ \alpha\in\real}. $
\begin{figure}[h!]
\centering
\includegraphics[width=0.475\linewidth]{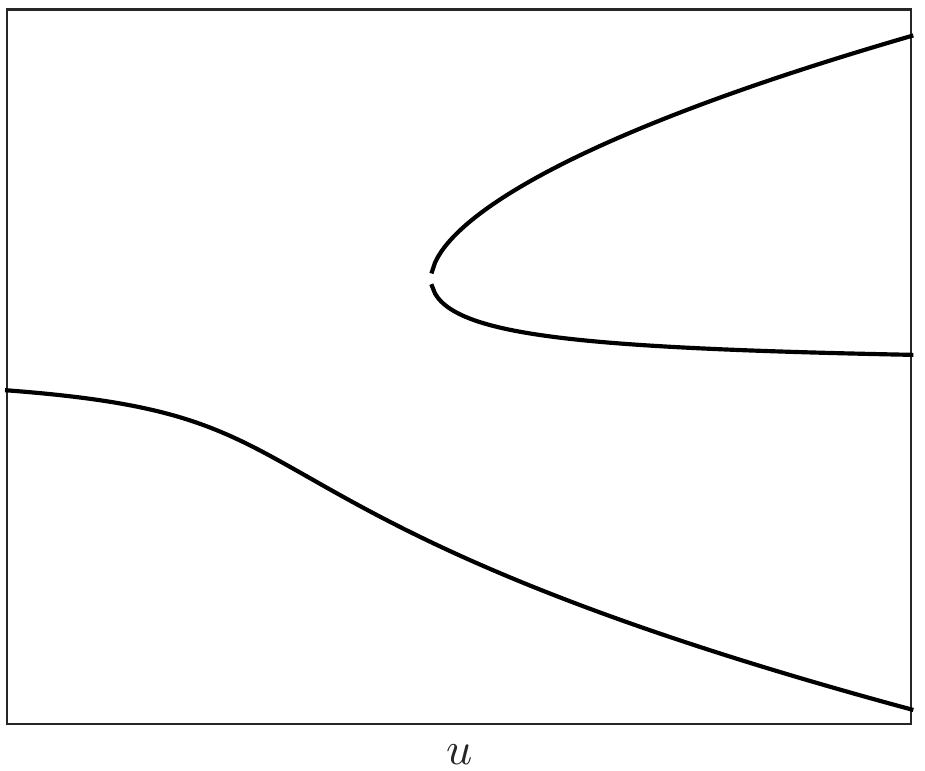} 
\includegraphics[width=0.475\linewidth]{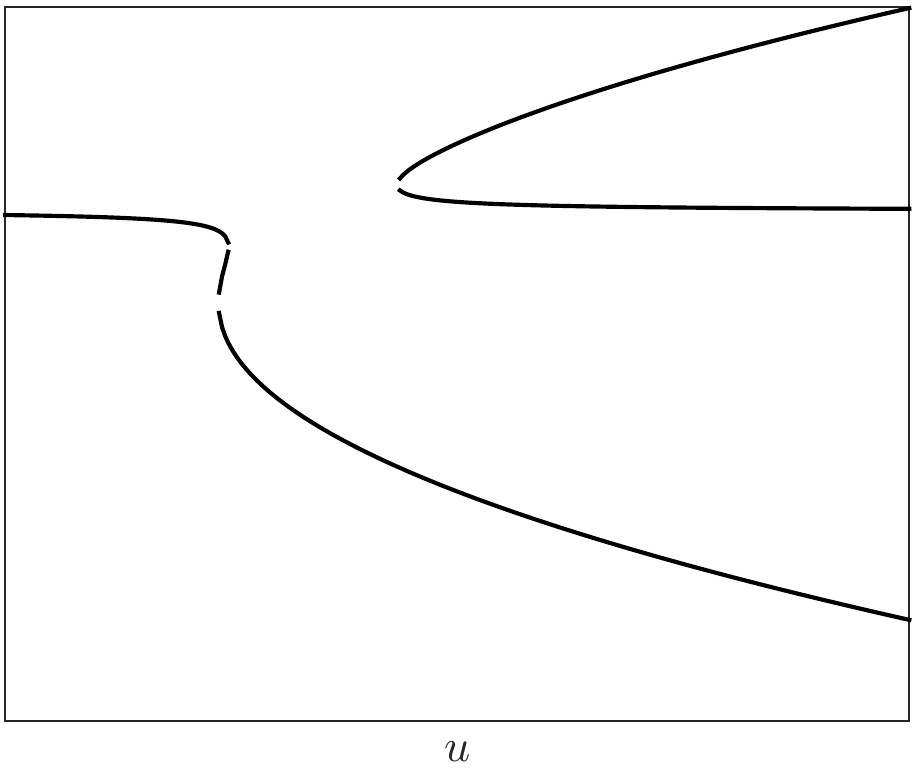}
\includegraphics[width=0.475\linewidth]{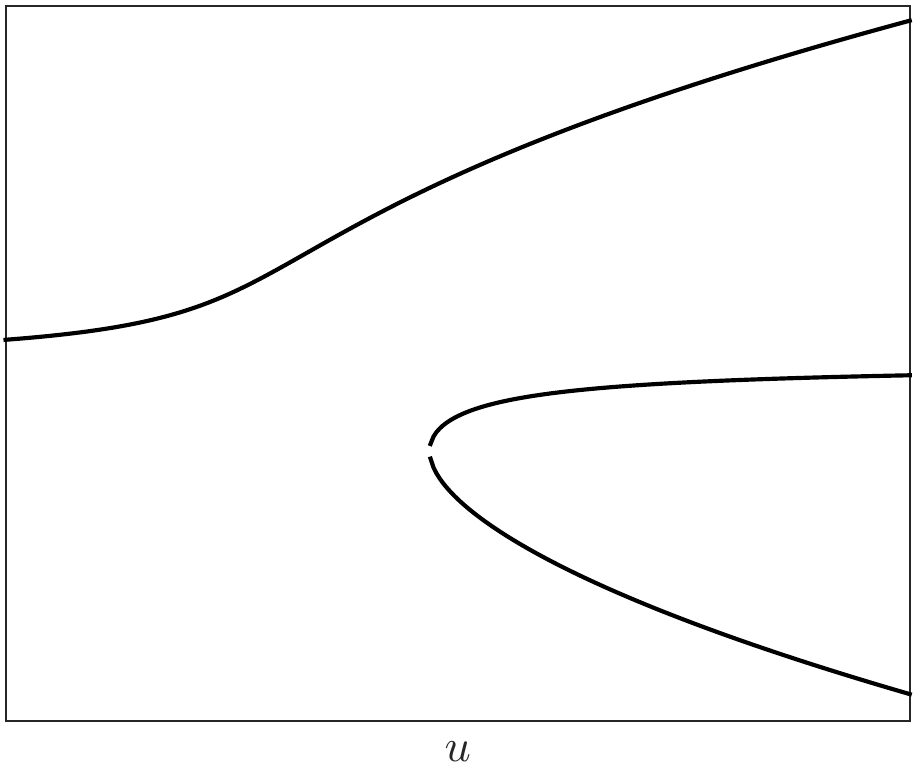}
\includegraphics[width=0.475\linewidth]{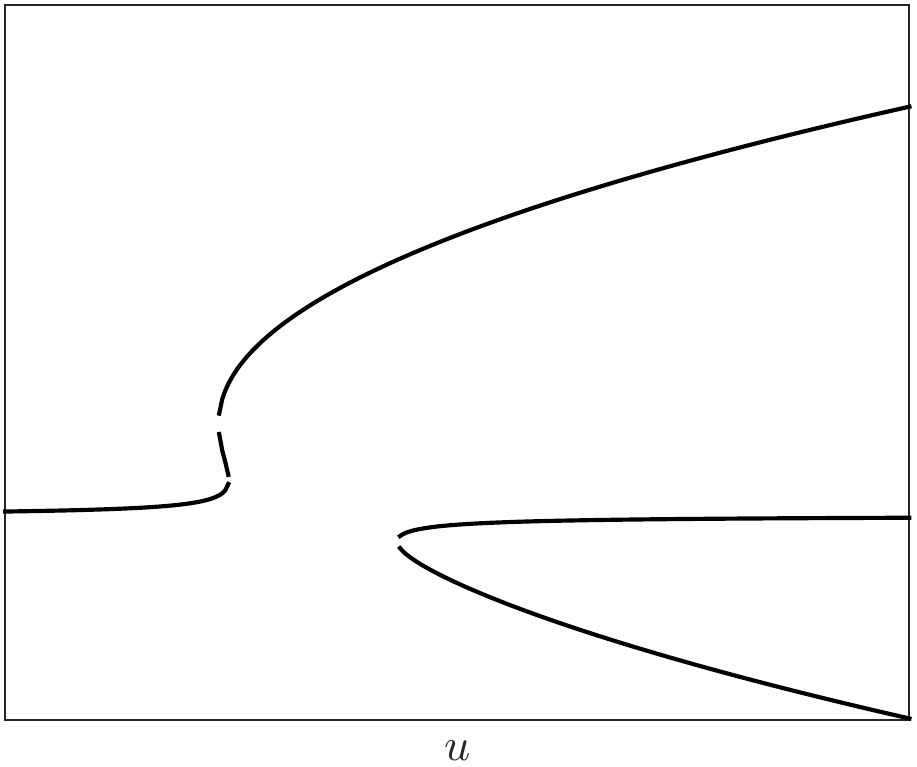}
\caption{Universal unfolding of the pitchfork: persistent bifurcation diagrams}\label{fig:univ-unfolding}
\end{figure}

If the bifurcation problem $g$ is recognized as the pitchfork bifurcation, then its universal unfolding $G$ is strongly equivalent to the $2$-parameter family of bifurcation problems $y^3 - u y + \alpha_1 + \alpha_2 y^2$, $\alpha_1, \alpha_2 \in \real$. Depending on the value of the unfolding parameters $\alpha_1$ and $\alpha_2$, the associated bifurcation diagram is (modulo strong equivalence) one of the persistent bifurcation diagrams of the pitchfork sketched in Fig.~\ref{fig:univ-unfolding}.

\subsection{Lyapunov-Schmidt reduction}

Consider the dynamical system 
\begin{equation}\label{eq:general-system}
\dot x_i = F_i (\bs x, u, \bs \alpha), \;\;  i=1, \ldots, N. 
\end{equation}
Let $(\bs x^*, u^*, \bs \alpha^*)$ be a fixed point~\eqref{eq:general-system}, and let $\mc L = \frac{\partial \bs F}{\partial \bs x}$, where $\bs F$ is the vector of $F_i,\; i=1, \ldots, N$. Suppose ${\rm rank }\;\mc L = n-1$. Let $E$ be the projection operator from $\real^n$ onto ${\rm Im }\; \mc L$. The basic idea of the Lyapunov-Schmidt reduction is to write the fixed points equation of~\eqref{eq:general-system} into the following two complementary equations
\begin{subequations}
\begin{align}
E \bs F (\bs x,  u,\bs \alpha) &=0, \label{eq:ls-1}\\
(I-E) \bs F (\bs x, u, \bs \alpha) &=0,\label{eq:ls-2}
\end{align}
on ${\rm Im }\; \mc L$ and ${\rm ker } \; \mc L$, respectively, and to write $\bs x = \bs v + \bs w$, with $\bs v \in {\rm ker } \; \mc L$ and $\bs w \in ({\rm ker }\; \mc L)^\perp$. 
\end{subequations}
Because for $|u-u^*|$ small $E\mc L$ is of full rank on ${\rm Im }\; \mc L$, by using the implicit function theorem~\eqref{eq:ls-1} can be solved to obtain $\bs w = \bs W(\bs v, u,  \bs \alpha)$. Substituting $\bs w$ in~\eqref{eq:ls-2}, we get 
\[
\bs f(\bs v, u, \bs \alpha) = (I - E) \bs F (\bs v +\bs W(\bs v, u, \bs \alpha), u, \bs \alpha),
\]
which maps to a one-dimensional subspace of $\real^n$ and its zeros are in one-to-one correspondence with the fixed points of~\eqref{eq:general-system}.

Let $\bs v_0$ and $\bar{\bs v}_0$ be vectors in ${\rm ker } \; \mc L$ and $({\rm Im }\;  \mc L)^\perp$. The vector equation $\bs f=0$ can be reduced to a scalar equation $G=0$ defined by
\[
G (y, u, \bs \alpha) = \bar{\bs v}_0 ^T \bs f(y \bs v_0, u,  \bs \alpha). 
\] 
Again, the zeros of $G$ are in one-to-one correspondence with the fixed points of~\eqref{eq:general-system}.  

The Lyapunov-Schmidt reduction transforms the bifurcation analysis of~\eqref{eq:general-system} in a neighborhood of $(x^*, u^*, \bs \alpha^*)$ to the scalar bifurcation problem $G (y, u, \bs \alpha)=0$. The tools from the previous section can then be used to ascertain the presence and nature of bifurcations.

\section{Analysis of the uninformed opinion dynamics}
\label{SEC: gloabl local analysis}

In this section, we analyze system~\eqref{EQ: LTIs network}. We summarize the stability properties of the equilibrium points of~\eqref{EQ: LTIs network} in the following theorem. For $u>1$,  let $\{0, \pm y^s\}$ be the three roots of the equation $-y+uS(y)=0$.  

\begin{theorem}[\bit{Uninformed opinion dynamics}]\label{thm: main unp pitch}
For the uninformed opinion dynamics~\eqref{EQ: LTIs network}, the following statements hold:
\begin{enumerate}
\item For $0<u< 1$, the origin is globally asymptotically stable and locally exponentially stable; 
\item At $u=1$, the origin is globally asymptotically stable and undergoes a pitchfork bifurcation; 
\item For $u>1$ and $|u-1|$ sufficiently small, there are exactly three fixed points at origin and $\pm y^s \bs 1_N$. The origin is a saddle point with a one dimensional unstable manifold. The stable manifold of the origin separates $\mathbb R^N$ into the basins of attraction of two exponentially stable symmetric stationary states $\pm y^s \bs1_N$. In particular, almost all (in topological and measure sense) trajectories converges to $\{y^s \bs 1_N\}\cup\{-y^s \bs 1_N\}$
\end{enumerate}
\end{theorem}

\begin{figure}[h!]
\includegraphics[width=0.45\textwidth]{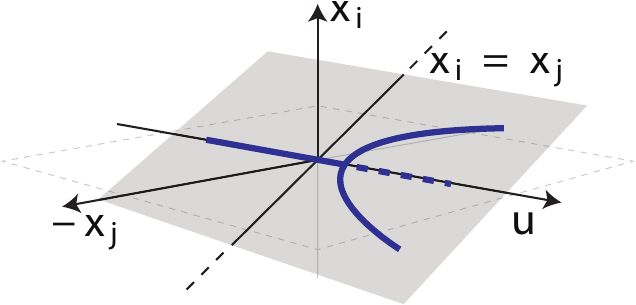}
\caption{For $u=1$, uninformed opinion dynamics~(\ref{EQ: LTIs network}) possess a pitchfork singularity at the origin. The steady state branches emerging at this singularity lie on the consensus space $\{x_i=x_j$, for all $i,j=1,\ldots,N\}$. Branches of the stable stationary states are depicted as solid lines. Branches of unstable stationary states are depicted as dashed lines. The consensus manifold is depicted as a grey plane.}
\label{FIG: pitch on consensus}
\end{figure}

Theorem \ref{thm: main unp pitch} is illustrated in Fig.~\ref{FIG: pitch on consensus}. The uninformed opinion dynamics~(\ref{EQ: LTIs network}) are symmetric, in the sense that the agents have no preferences toward the positive or negative alternative. At a mathematical level, this is reflected in the odd symmetry of the vector field and in the fact that the pitchfork is symmetric, that is, the two stable steady states emerging at the pitchfork are mirror-symmetric with respect to the origin and belong to the consensus manifold. 

We devote the remainder of the section to prove Theorem~\ref{thm: main unp pitch}. We establish the theorem for the case of $0 < u \le 1$ using Lyapunov analysis. The case with $u>1$ is more involved. In this case, we resort to singularity analysis to show that dynamics~\eqref{EQ: LTIs network} exhibit a pitchfork bifurcation at $u=1$ and have only three equilibrium points for $u-1$ sufficiently small. Finally, we use a result from monotone systems theory to establish almost global convergence to these equilibrium points. 


\subsection{Analysis for the case $0<u \le 1$}
Consider the quadratic Lyapunov function $V(\bs x)=\frac{1}{2}\bs x^T \bs x$:
\begin{IEEEeqnarray}{rCl}\label{EQ: lyap}
\IEEEyesnumber
\dot V&=&\bs x^T(-D\bs x+u A \bs S(\bs x))\nonumber\\
&=&\bs x^T(-D \bs x+uD \bs S(\bs x)-u D \bs S(\bs x)+u A \bs S(\bs x))\nonumber\\
&=&- \bs x^TD(\bs x-u \bs S(\bs x))- u \bs x^TL \bs S(\bs x))\nonumber\\
&<&-u \bs x^TL \bs S(\bs x)\nonumber\\
&\leq&0,\quad \forall \bs x\neq0,
\end{IEEEeqnarray}
where we have used the fact that $uS$ is a monotone function in the sector $[0,1]$, that $D$ is a diagonal Hurwitz matrix, and that $L$ is a positive semidefinite matrix. 

Local exponential stability for $0<u<1$ follows by noticing that the linearization of (\ref{EQ: LTIs network}) at the origin is
$$\dot{\delta \bs x}=(-D+uA)\delta \bs x, $$
and $(-D+uA)$ is a diagonally dominant and Hurwitz matrix.

%

\subsection{Analysis of the case $u>1$}

\subsubsection*{\bit{Lyapunov Schmidt reduction}}

We start by applying Lyapunov-Schimdt reduction to the fixed point equation of~\eqref{EQ: LTIs network}
\begin{equation}\label{EQ: LTIs network fp}
-D \bs x+u A \bs S(\bs x)=0,
\end{equation}
near $(x,u)=(0,1)$.

For $u=1$, the linearization of (\ref{EQ: LTIs network fp}) is $-L\bs x=0$, with $L=-A + D$ the network Laplacian matrix and, by the strongly connected assumption, ${\rm ker}\; L={\rm span}\; \bs 1_N $. Let $\Pi=I_N-\frac{1}{N} \bs 1_N \bs 1_N^T$ be the projector from $\real^N$ onto $\bs 1_N^\perp$ and define $\bs w=\Pi  \bs x$ and $(I-\Pi)\bs x=: y \bs 1_N$, for some $y\in\mathbb R$, the projections of $ \bs x$ on $\bs 1_N^\perp$ and $\bs 1_N$, respectively.  We split (\ref{EQ: LTIs network fp}) into two reduced equations
\begin{IEEEeqnarray}{rCl}\label{EQ: LTIs network fp reduced 1}
\IEEEyesnumber
\Pi(-D(\bs w+y \bs 1_N)+uA \bs S(\bs w+y \bs 1_N))&=&0,\IEEEyessubnumber\\ 
(I-\Pi)(-D(\bs w+y \bs 1_N)+uA \bs S(\bs w+y \bs 1_N))&=&0,\IEEEyessubnumber \label{EQ: LTIs network fp reduced 1.1}
\end{IEEEeqnarray}
on $\bs 1_N^\perp$ and ${\rm span}\; \bs 1_N$ respectively. Since the linearization of (\ref{EQ: LTIs network fp}) restricted to $\bs 1_N^\perp$ is non-singular for $u$ sufficiently close to 1, by the implicit function theorem we can solve (\ref{EQ: LTIs network fp reduced 1}a) for $\bs w$ locally around the origin as follows.

Using the Taylor series expansion in $\bs x$ and using $L \bs 1_N=0$, we can rewrite (\ref{EQ: LTIs network fp reduced 1}a) as
\begin{equation*}
(-\Pi L+(u-1)\Pi A)\bs w+y(u-1)\Pi A \bs 1_N + \mathcal O(u \bs x^3)=0,
\end{equation*}
where $\mc O(u \bs x^3)$ denotes the vector with third and higher order terms in the Taylor series expansion. 

For $u$ sufficiently close to $1$, the matrix $\Pi L-(u-1)\Pi A$ is invertible on $\bs 1_N^\perp$. Therefore, the first order solution is $\bs w=(u-1)y \bs c$ with
\begin{equation*}
\bs c= (\Pi L-(u-1)\Pi A)^{-1}\Pi A \bs 1_N.
\end{equation*}
Invoking the fact that (\ref{EQ: LTIs network fp reduced 1}a) is odd, the next correction term is of the third order. Thus, the solution to (\ref{EQ: LTIs network fp reduced 1}a) is given by
\begin{equation}\label{EQ: LTIs network fp reduced 2}
\bs w=(u-1)y \bs c+\tilde{\bs w},\ \ \bs c,\tilde{\bs w} \in \bs 1_N^\perp,
\end{equation}
with $\tilde w_i=\mathcal O(uy^3)$ and smooth for all $i=1,\ldots,N$. 

We now substitute (\ref{EQ: LTIs network fp reduced 2}) into (\ref{EQ: LTIs network fp reduced 1}b) and obtain
\begin{multline}\label{EQ: vector fp eq unpert}
(I-\Pi)\Big(-D\big(y \bs 1_N+(u-1)y \bs c+\tilde {\bs w}\big) \\
+ uA \bs S\big(y \bs 1_N+(u-1)y\bs c+\tilde {\bs w} \big)\Big) = 0.
\end{multline}
Since $(I- \Pi)$ is a projector onto $\bs 1_N$, each entry in the column vector in the LHS of~\eqref{EQ: vector fp eq unpert} is identical and equal to 
\begin{multline*}
-\frac{1}{N}\sum_{i=1}^N \Big(d_iy + (u-1)y d_ic_i + d_i\tilde w_i+\\
 -\sum_{j=1}^N ua_{ij}  S\big(y+(u-1)yc_j+\tilde w_j\big)\Big).
\end{multline*}
Recalling that the interconnection graph is balanced, and performing some algebraic simplifications, reduces the above expression to
\begin{multline}\label{EQ: scalar fp eq unpert}
g(y,u) := -\frac{1}{N}\sum_{i=1}^N d_i \Big(y+(u-1)yc_i+\tilde w_i  \\
- uS(y+(u-1)yc_i+\tilde w_i)\Big).
\end{multline}

\subsubsection*{\bit{Singularity Analysis}}
We now perform singularity analysis and solve the recognition problem for the bifurcation problem $g$. 
To recognize a pitchfork bifurcation, we need to show that $g=g_y=g_{yy}=g_u=0$ and $g_{yyy}<0,\ g_{yu}>0$ at $(y,u)=(0,1)$. Clearly $g(0,1)=0$. 
For brevity,  we let $\varsigma =y+(u-1)yc_i+\tilde w_i$ and we omit the factors $-\frac{1}{N}\sum_{i=1}^Nd_i$. We have
\begin{align*}
g_y &=1+(u-1)c_i+\frac{\partial\tilde w_i}{\partial y}
-u S'(\varsigma)\Big(1+(u-1)c_i+\frac{\partial\tilde w_i}{\partial y}\Big), \\
\!g_{yy}&= \!\frac{\partial^2\tilde w_i}{\partial y^2}\! - \!uS''(\varsigma)\Big(\!1+u(u-1)c_i+\frac{\partial\tilde w_i}{\partial y}\! \Big)^2 \!\! - \! uS'(\varsigma)\frac{\partial^2\tilde w_i}{\partial y^2}, \\
\!\!g_{yyy} \! &=\frac{\partial^3\tilde w_i}{\partial y^3}-uS'''(\varsigma)\Big(1+(u-1)c_i+u\frac{\partial\tilde w_i}{\partial y}\Big)^3\\
& \quad -uS'(\varsigma)\frac{\partial^3\tilde w_i}{\partial y^3} -uS''(\varsigma)2\Big(\!1+(u-1)c_i+\frac{\partial\tilde w_i}{\partial y}\!\Big)\frac{\partial^2\tilde w_i}{\partial y^2}\\
& \qquad \qquad \quad -uS''(\varsigma)\left(1+(u-1)c_i+\frac{\partial\tilde w_i}{\partial y}\right)\frac{\partial^2\tilde w_i}{\partial y^2} ,\\
g_u&=yc_i-S(\varsigma)-u S'(\varsigma)\Big(yc_i+\frac{\partial\tilde w_i}{\partial u}\Big), \quad \text{and}\\
g_{yu} & =c_i\! -\! S'(\varsigma)\Big(\!1+(u-1)c_i+\frac{\partial\tilde w_i}{\partial y}\!\Big) 
\!- \!uS'(\varsigma)\Big(\!c_i+\frac{\partial^2w_i}{\partial y\partial u}\!\Big) \\
&\quad -uS''(\varsigma)\Big(1+(u-1)c_i+\frac{\partial\tilde w_i}{\partial y}\Big)\Big(yc_i+\frac{\partial\tilde w_i}{\partial u}\Big).
\end{align*}
Recalling that $S'(0)=1$, $S''(0)=0$, $S'''(0)<0$, and since $\tilde w_i=\mathcal O(uy^3)$, $\frac{\partial\tilde w_i}{\partial y}=\frac{\partial^2\tilde w_i}{\partial y\partial u}=0$, we have at $(y,u)=(0,1)$:
\begin{align*}
g &=g_y=g_{yy}=g_u=0, \quad \text{and}\\
g_{yyy}&=\frac{1}{N}\sum_{i=1}^Nd_iS'''(0)<0,\quad g_{uy}=\frac{1}{N}\sum_{i=1}^Nd_i>0.
\end{align*}
This completes the recognition of the pitchfork bifurcation. 


%
%

\subsubsection*{\bit{Equilibrium points for $u>1$}}

The dynamical system~\eqref{EQ: LTIs network} exhibits a pitchfork bifurcation at $u=1$. Hence, by the Lyapunov-Schmidt reduction and singularity analysis, for $u>1$ and $u-1$ sufficiently small,~\eqref{EQ: LTIs network} has at exactly three equilibrium points in a sufficiently small neighborhood of the origin. One of these equilibrium points is at origin.  By replacing $\bs x=\pm y^s \bs 1_N$, the $i$-th line of equation (\ref{EQ: LTIs network fp}) reads
\begin{IEEEeqnarray*}{rCl}
-d_i\pm y^s+\sum_{j=1}^Nua_{ij}S(\pm y^s)&=&d_i(\mp y^s+uS(\pm y^s)) =0, 
\end{IEEEeqnarray*}
which proves that $x^s_\pm=\pm y^s1_N$ are the other fixed points of (\ref{EQ: LTIs network}). By \cite[Theorem I.4.1]{Golubitsky1985}, the fixed point at the origin has exactly one positive eigenvalue (corresponding to the eigenvalue crossing zero at the pitchfork) and, therefore, it is a saddle with a one-dimensional unstable manifold. The other two fixed points $\pm y^s \bs 1_N$ have the same local stability properties of the bifurcating fixed point, that is, they are locally exponentially stable.

We now show that these three equilibrium points are (globally) the only fixed points for $u-1$ sufficiently small. 
Let $0\in B\subset\mathbb R^N$ be a (sufficiently small) neighborhood of the origin where thesingularity analysis hold. Since for $u=1$ there are no fixed points on $\mathbb R^N\setminus B$; by continuity of the vector field on $u$ and because $|u-1|\bs S(\bs x)$ is bounded on $\mathbb R^N$ the same must be true for $|u-1|$ sufficiently small. Since the only fixed points inside $B$ are those predicted by singularity theory, our claim is proved.

\subsubsection*{\bit{Almost global convergence}}

We now show that the system~\eqref{EQ: LTIs network} is a strongly monotone system, and then use a result from monotone system theory that ensures that almost all trajectories of~(\ref{EQ: LTIs network}) converge to the set of these three equilibrium points.

For $u>0$, uninformed opinion dynamics~(\ref{EQ: LTIs network}) is a {\it cooperative}\footnote{The exact definition of \emph{cooperation} is a debated topic in the ecology literature. We use the term \emph{cooperative} in a strictly mathematical sense as defined in~\cite{Hirsch1988}.} system~\cite[Page 1]{Hirsch1988}, that is, its Jacobian matrix satisifies the condition
$$\frac{\partial\dot x_i}{\partial x_j}(\bs x)=a_{ij}uS'(x_j)\geq 0,\; \forall i\neq j, \; \forall \bs x\in\mathbb R^N.$$
Moreover, the Jacobian matrix of (\ref{EQ: LTIs network}) is also irreducible\footnote{A square matrix is irreducible if and only if  the associated directed graph is strongly connected.} for each $\bs x\in\mathbb R^N$. Indeed, because $uS'(x_j)>0$ for all $j=1,\ldots,N$ and all $\bs x\in\mathbb R^N$, the graph associated \cite[Definition 6.2.11]{Horn1985} with the Jacobian matrix of (\ref{EQ: LTIs network}) is the same as the graph associated with $A$ for all $\bs x\in\mathbb R^N$. Because this graph is strongly connected, the Jacobian matrix of (\ref{EQ: LTIs network}) is irreducible for all $\bs x\in\mathbb R^N$ \cite[Theorem 6.2.24]{Horn1985}.

The flow of a cooperative system with an irreducible Jacobian matrix is {\it strongly monotone} \cite{Hirsch1988}. In particular, the flow induced by (\ref{EQ: LTIs network}) is strongly monotone on the whole $\mathbb R^N$. Moreover, because $-D$ is Hurwitz and $uA\bs S(\bs x)$ is bounded, (\ref{EQ: LTIs network}) has no unbounded trajectories. Without recalling technicalities and definitions about monotone systems (the reader is referred to \cite{Hirsch1988}), we now specialize \cite[Theorem 0.1]{Hirsch1988} to model (\ref{EQ: LTIs network}).

For a given $u \in \real_{>0}$, let $\mathcal S\subset\mathbb R^N$ denote the set of fixed points of model (\ref{EQ: LTIs network}). Let $\mathcal C\subset\mathbb R^N$ denote the set of points whose $\omega$-limit set is in $\mathcal S$ ($\bs x\in\mathcal C\ \Leftrightarrow\ \omega(\bs x)\subset\mathcal S$). Then, 
\begin{enumerate}
\item $\mathcal C$ is residual in $\mathbb R^N$, i.e., its complement is nowhere dense;
\item $\mu(\mathbb R^N\setminus\mathcal C)=0$ for all Gaussian measures $\mu$.
\end{enumerate}
Now, the almost global convergence to $\{\pm y^s \bs 1_N\}$ follows immediately.

\section{Analysis of the informed opinion dynamics: Complete graph}
\label{sec:unfold}

In this section, we analyze the informed opinion dynamics~\eqref{eq:opinion-forced} for the complete graph. 
We leave the analysis  for a generic strongly connected balanced graph as future work. For a complete graph,~\eqref{eq:opinion-forced} reduces to
\begin{equation}\label{EQ: ata network with prefs}
\dot x_i=-(N-1)x_i+\sum_{\substack{j=1\\j\neq i}}^NuS(x_j)+\alpha_i,\quad i=1,\ldots,N.
\end{equation}

We summarize the main result of this section in the following proposition. 
\begin{proposition}[\bit{Informed opinion dynamics}] \label{prop:forced-opinion-dynamics}
The following statement holds for the informed opinion dynamics~\eqref{eq:opinion-forced} defined on a complete graph: 
the Lyapunov-Schmidt reduction of the fixed points equation of ~\eqref{eq:opinion-forced} around $(\bs x^*, u, \bs \alpha) =(0, 1, 0)$ is 
\begin{align}\label{EQ: ata scalar fp eq}
G(y,u; \bs \alpha)=-y+\langle \bs \alpha\rangle+\langle\tilde{\bs w} \rangle +\frac{1}{N}\sum_{i=1}^NuS(\alpha_i^\perp+y+\tilde w_i),
\end{align}
where $\bs \alpha^\perp:=\Pi \bs \alpha$, $\tilde w_i=\mathcal O(uy^3)$, and $\langle \bs z\rangle=\frac{1}{N}\sum_{j=1}^Nz_j$ for all $\bs z\in\mathbb R^N$.
\end{proposition}

{\bf Proof.}
We start by applying the Lyapunov-Schmidt reduction to the $N$-dimensional fixed point equation of~\eqref{EQ: ata network with prefs}
\begin{equation}\label{EQ: ata network fp}
-(N-1)\bs x+(\bs 1_N \bs 1_N^T- I_N)u \bs S(\bs x)+ \bs \alpha=0.
\end{equation}
Similar to the previous section, we define $\bs w=\Pi \bs x$ and $(I-\Pi)\bs x=:y \bs 1_N$, and we split (\ref{EQ: ata network fp}) into two equations
\begin{subequations}
\begin{multline}\label{EQ: ata network fp reduced 1}
\Pi\big(\!-\!(N-1)(\bs w+y\bs 1_N) \\
+ (\bs 1_N \bs 1_N^T-I_N)uS(w+y \bs 1_N)+\bs \alpha\big)=0,
\end{multline}
\begin{multline} \label{EQ: ata network fp reduced 1.1}
(I\!-\!\Pi)\big(\!-\!(N-1)(\bs w+y \bs 1_N)c \\
+(\bs 1_N \bs 1_N^T-I_N)uS(w+y \bs 1_N)+\bs \alpha\big)=0,
\end{multline}
\end{subequations}
defined on $\bs 1_N^\perp$ and $\bs 1_N$, respectively. Similar to the previous section, equation~(\ref{EQ: ata network fp reduced 1}) is non-singular on $\bs 1_N^\perp$ for $|u-1|$ and $|\alpha|$ sufficiently small.

Using the Taylor series expansion in $\bs x$, and recalling that $\Pi \bs 1_N=0$ and $\Pi \bs w= \bs w$, (\ref{EQ: ata network fp reduced 1}) can be rewritten as
\begin{align*}
&\Pi(\bs 1_N \bs 1_N^T-I_N)u(\bs w+y \bs 1_N)\!-\!(N\!-\!1) \bs w+\Pi \bs \alpha+\mathcal O((u \bs x)^3)\\
&=-(N+u-1) \bs w+\Pi \bs \alpha+\mathcal O((u \bs x)^3) = 0,
\end{align*}
where we used $\Pi(\bs 1_N \bs 1_N^T-I_N) = -\Pi$. It follows that, to the first order, the solution is $\bs w=\frac{\Pi \bs \alpha}{N+u-1}$ and consequently, 
\begin{equation}\label{EQ: ata network fp reduced 2}
\bs w=\frac{\Pi \bs \alpha}{N+u-1}+\tilde{\bs w},\quad \tilde{\bs w} \in \bs 1_N^\perp,
\end{equation}
with $\tilde w_i=\mathcal O(uy^3)$ and smooth for all $i=1,\ldots,N$.

%
%

Substituting (\ref{EQ: ata network fp reduced 2}) in (\ref{EQ: ata network fp reduced 1.1}) and proceeding as in the previous section, we obtain the scalar fixed point equation~\eqref{EQ: ata scalar fp eq}.   \oprocend


For $\bs \alpha=0$, (\ref{EQ: ata scalar fp eq}) is a special case of (\ref{EQ: scalar fp eq unpert}). Therefore, by Theorem \ref{thm: main unp pitch}, the bifurcation problem $G$ also exhibits a pitchfork bifurcation for $\bs \alpha =0$. A nonzero $\bs \alpha$ has the effect of {\it unfolding} the symmetric pitchfork.


We now briefly discuss the effect of distribution of $\bs \alpha$ across different nodes in the graph. 
It is easy to see using the expression in~\eqref{EQ: ata scalar fp eq} that to the first order the contribution of $\bs \alpha^\perp$ vanishes.  This suggests that the unfolding of the pitchfork in (\ref{EQ: ata scalar fp eq}) solely depends on the average of the preference vector $\langle\bs \alpha\rangle$. This reflects the symmetry of (\ref{EQ: ata network with prefs}) with respect to agent permutation. 

However, the fact that $\bs \alpha^\perp$ does not influence the bifurcation problem is only partially true. The components $\alpha^\perp_i$, while not unfolding the pitchfork, move the bifurcation point. That is, in the presence of a nonzero $\bs \alpha^\perp$, the pitchfork singularity occurs at $u=u^*>1$ even if $\langle\bs\alpha\rangle=0$. 

For $\langle \bs \alpha\rangle=0$, evaluating the derivative of (\ref{EQ: ata scalar fp eq}) with respect to $y$ at $y=0$ gives
\begin{IEEEeqnarray*}{rCl}
G_y(0,u;\alpha)&=&-(N-1)+\frac{N-1}{N}\sum_{i=1}^NuS'(\alpha_i^\perp),
\end{IEEEeqnarray*}
where we used that $\tilde w_i=\frac{\partial w_i}{\partial y}=0$. 

For $\bs \alpha^\perp=0$, we get, as expected, $G_y(0,1;0)=0$. However, because $S'(z)<1$ for $z\neq0$, we immediately see that, for $\bs \alpha^\perp\neq 0$, $u>1$ is necessary to enforce the singularity condition $G_y=0$. The following example verifies this prediction.

\begin{example}
For $N=5$, $\bs \alpha=[1\ \ 1\ -1\ -1\ \ 0]$, and $S(\cdot)=\tanh(\cdot)$, numerical simulations shows that the pitchfork singularity exists for $u^*>1.034$ \oprocend
\end{example}

The development above is not rigorous, mainly because for $\bs \alpha^\perp\neq0$ the bifurcation point is not at the origin. However, it allows us to formulate the following claim, which will be proved in an extended version of this work.

\begin{claim}\label{CLM: ata unfolding structure}
There exists $\delta_1>0$, such that, if $0< \|\bs \alpha^\perp\|<\delta_1$ and $\langle \bs \alpha\rangle=0$, then there exist $u^*>1$ and $\bs x^*\in\mathbb R^N$, such that the informed opinion dynamics~(\ref{EQ: ata network with prefs}) exhibits a pitchfork singularity at $(\bs x,u)=(\bs x^*,u^*)$. In particular, the Lyapunov-Schmidt reduction of its fixed point equation (\ref{EQ: ata network fp}) at $(\bs x,u)=(\bs x^*,u^*)$ is strongly equivalent to the pitchfork singularity.
\end{claim}

We now illustrate the prediction of Claim~\ref{CLM: ata unfolding structure} using an extended example in the next section. 

\section{Collective Decision-Making in Honeybees}
\label{SEC: honey bee}

In the honeybee nest site selection problem, an entire swarm must unanimously choose a nest site where it will live as a new colony with its queen. The process starts with a number of bees each scouting out a possible nest site and ends with the swarm choosing the best of the scouted-out alternatives. Each informed honeybee scout uses explicit signaling in the form of a ``waggle dance" to recruit uninformed bees to commit to its discovered nest site. Seeley et al.~\cite{Seeley12}  conducted experiments in which honey-bees had to choose between two alternative nest sites. Using empirical data, they showed that, along with the waggle dance,  the scouts also use a cross-inhibitory stop-signal to
make the scouts recruiting for the competing alternative stop dancing. This stop-signal contributes positively to the collective decision-making; of particular note, it was shown to help with breaking deadlock in the case of two near equal value alternative sites.

Consider the following  informed-opinion dynamics
\begin{align}\label{eq:agent-model-input}
\dot{\bs x}(t) = - u_1 D \bs x(t) +  u_2 A S( \bs x) + \bs \alpha,
\end{align}
where $\bs \alpha \in \mathbb R^n$ is the vector with entries equal to the value of the nest site that the corresponding agent is exposed to. We associate positive (negative) values with the first (second) nest site and the value absolute value with the nest quality.  The entries of $\bs \alpha$ associated with the uninformed agents are zero. 

In~\eqref{eq:agent-model-input}, $u_1\in \real_{>0}$ is the rate at which an agent commits to the nest site it is exposed to, and $u_2\in \real_{>0}$ is the rate at which an agent is attracted  by other agents to the nest sites. Notice that ~\eqref{eq:agent-model-input} is equivalent to the informed opinion dynamics~\eqref{eq:opinion-forced} defined on a  time-scale $u_1 t$ and with $u=u_2/u_1$.  Here, we focus on complete graphs and in this case~\eqref{eq:agent-model-input} is equivalent to~\eqref{EQ: ata network with prefs}.

Suppose honeybees have to select between two nest sites with values $v_1\in (0,1)$ and $v_2 \in (-1,0)$, respectively. Then, $\alpha_i=v_1$ if the bee $i$ is exposed to the nest with value $v_1$, $\alpha_i=v_2$ if the bee $i$ is exposed to the nest with value $v_2$, and $\alpha_i=0$ if the bee $i$ is not exposed to any nest site. Further assume that the number of agents exposed to each nest site is the same. Let $\bar v = (v_1+v_2)/2$ and $u_1 = 1/ \bar v$. Such choice of $u_1$ means that for small values of $\bar v$ and a fixed value of $u_2$, the opinions converge to zero quickly. 

\begin{figure}[h!]
\centering
\subfigure[Opinion of Population]{\includegraphics[width=0.49\linewidth]{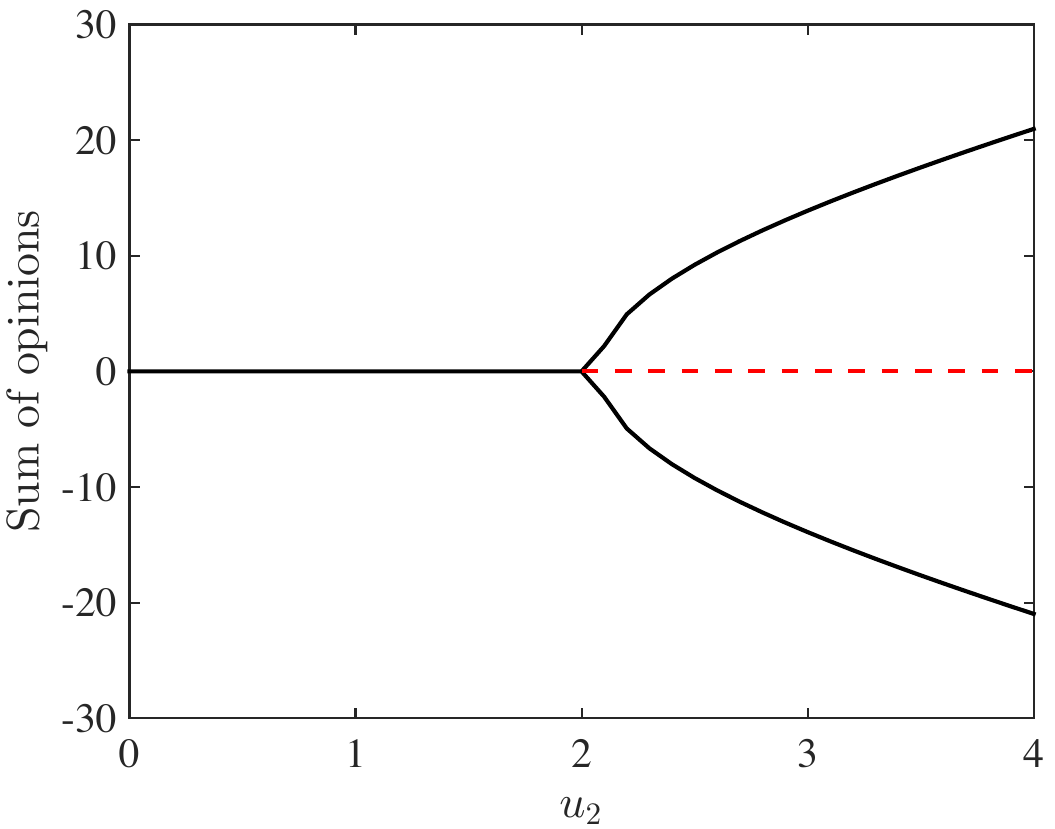}\label{fig:pop-opinion}}
\subfigure[Bifurcation Point]{\includegraphics[width=0.49\linewidth]{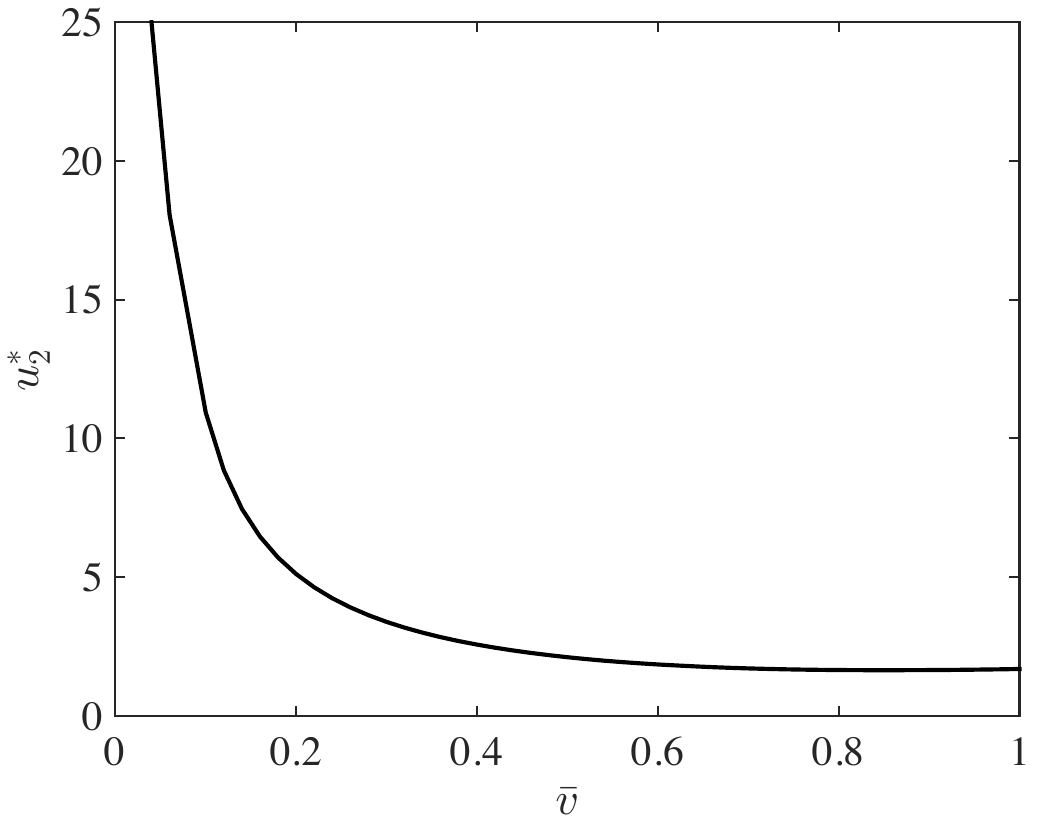}\label{fig:bif-pt-val}}
\subfigure[Unfolding]{\includegraphics[width=0.49\linewidth]{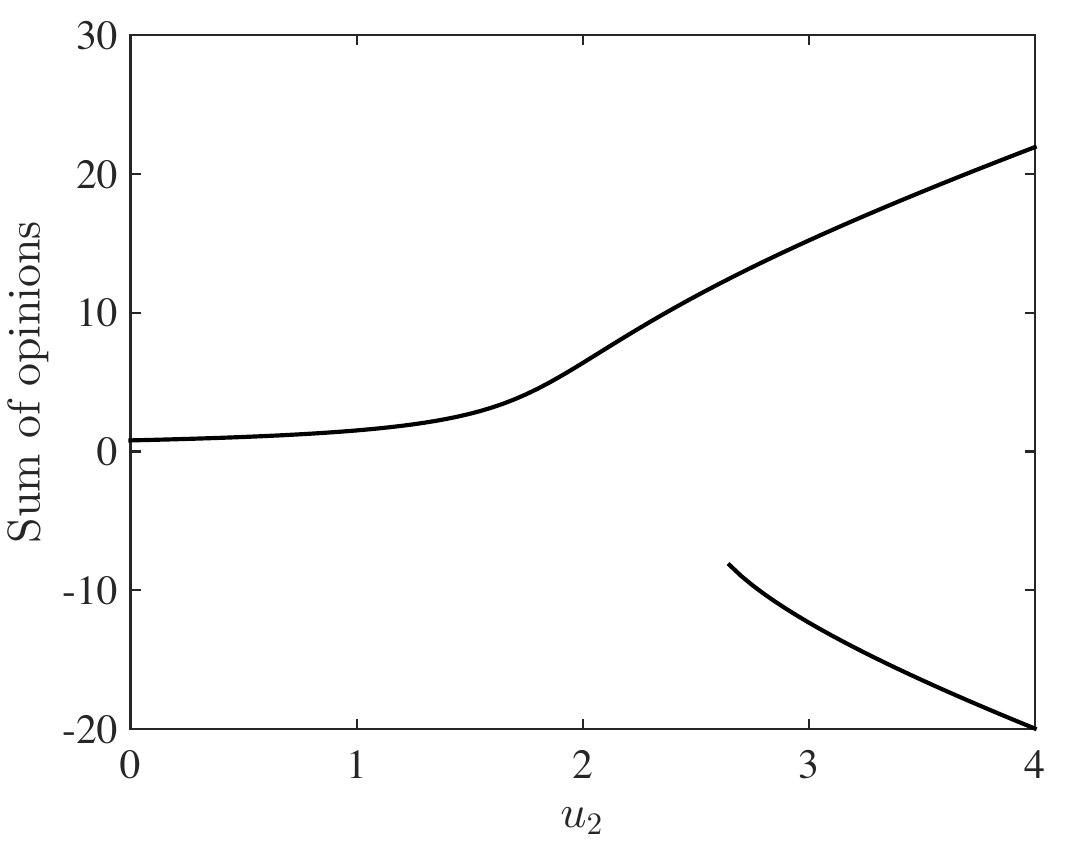} \label{fig:unfolding}}
\caption{Decision-making in selecting nest sites. (a) Bifurcation diagram  for an all-to-all graph with $11$ agents and value of nest site equal to $0.5$. Both nest sites were exposed to four agents each. (b) Bifurcation point as a function of the value of nest sites. (c) Unfolding of the pitchfork bifurcation for asymmetric nest sites. Only stable equilibrium points are shown. The value of nest sites are $v_1=0.7$ and $v_2=0.5$.  \label{fig:symm}}
\end{figure}

We first consider the symmetric case in which $v_1=v_2$. Following the discussion in \S\ref{sec:unfold}, for small values of $v_1$, the unfolding of the pitchfork depends on the average of $\bs \alpha$ which is equal to zero. Thus, in the symmetric case, there is no unfolding of the pitchfork. In this case for small values of $u_2$, the sum of all opinions converges to zero and  there is no unanimity in the group. However, at a critical value of $u_2$ a pitchfork bifurcation happens and the sum of opinion converges to a positive or a negative equilibrium point that corresponds to a  decision for one of the nest sites. The parameter $u_2$ maps to the cross-inhibition parameter in honeybees. The sum of opinions at steady state as a function of $u_2$ is shown in Fig.~\ref{fig:pop-opinion} and the bifurcation point as a function of the value of the nest sites is shown in Fig.~\ref{fig:bif-pt-val}. The variation of the bifurcation point as a function of the value of the nest site is qualitatively identical to the trend observed in the study of the honeybee model in~\cite{Pais2013}.

In the context of asymmetric nest site values the average of $\bs \alpha$ is not zero and the associated bifurcation diagram is shown in Fig.~\ref{fig:unfolding}. In this case for small cross-inhibition $u_2$, the sum of opinion of agents is zero and there is no unanimous decision. As we increase the value of $u_2$, an unfolding of pitchfork bifurcation is observed that favors the better nest site. A further increase in $u_2$ leads to emergence of another stable equilibrium point and in this case a decision can be made in favor of inferior nest site as well. Again this behavior is qualitatively identical to the behavior predicted by population-level models in~\cite{Pais2013}~and~\cite{Seeley12}.

\section{Conclusions and perspectives}
\label{SEC: conclusions}
In this paper, we developed agent-level mechanistic models for the realization of collective decision-making behavior in animal groups. We rigorously characterized the steady-state behavior of these models. We applied these models to the honeybee nest site selection problem and showed that our model captures the empirically observed behavior as well as the behavior predicted by  the population-level models. 

There are several possible extensions of this work. We analyzed the informed opinion dynamics model for the complete graph, and it is of interest to analyze the model for a generic strongly connected graph. 
In particular, the effect of graph topology of the onset of the bifurcation needs to be understood. Furthermore, the influence of the location of an informed agent in the network on the decision-making performance needs to be rigorously characterized. Another important future direction is to apply this model to other decision-making scenarios, e.g., food site selection in collective foraging and decision making in biological neuronal networks. 

\footnotesize
\bibliographystyle{unsrt}
\bibliography{FSL}

\end{document}